\documentclass[a4paper,11pt]{amsart}
\usepackage{amssymb,amsxtra,eucal}
\usepackage{fullpage}
\usepackage[all]{xy}

\theoremstyle{plain}
\newtheorem{thm}{Theorem}[section]

\newtheorem{prop}[thm]{Proposition}
\newtheorem{cor}[thm]{Corollary}

\theoremstyle{definition}
\newtheorem{rem}[thm]{Remark}
\newtheorem{defn}[thm]{Definition}

\newtheorem{assump}[thm]{Assumption}

\numberwithin{equation}{section}

\def\ie{\emph{i.e.}}

\def\:{\colon}
\def\.{\cdot}

\def\<{\left\langle}
\def\>{\right\rangle}
\def\({\left(}
\def\){\right)}

\def\epsilon{\varepsilon}
\def\phi{\varphi}

\def\leq{\leqslant}
\def\geq{\geqslant}

\def\lra{\longrightarrow}

\def\tilde#1{\widetilde{#1}}

\def\F{\mathbb{F}}

\def\Q{\mathbb{Q}}

\def\Z{\mathbb{Z}}

\def\Times_#1{\mathop{\times}_{#1}}

\def\oTimes_#1{\mathop{\otimes}_{#1}}
\def\oPlus_#1{\mathop{\bigoplus}_{#1}}

\DeclareMathOperator{\Hom}{Hom}

\def\BP{\mathit{BP}}
\def\MU{\mathit{MU}}
\def\bpn{\BP\langle n\rangle}
\def\bp1{\BP\langle 1\rangle}

\begin{document}
\title[Uniqueness of $E_\infty$ structures]
{Uniqueness of $E_\infty$ structures for connective covers}
\author{Andrew Baker \and Birgit Richter}
\address{Mathematics Department, University of Glasgow,
Glasgow G12 8QW, Scotland.}
\email{a.baker@maths.gla.ac.uk}
\urladdr{http://www.maths.gla.ac.uk/$\sim$ajb}
\address{Fachbereich Mathematik der Universit\"at Hamburg,
20146 Hamburg, Germany.}
\email{richter@math.uni-hamburg.de}
\urladdr{http://www.math.uni-hamburg.de/home/richter/}
\thanks{
We are grateful to John Rognes who suggested to exploit the
functoriality of the connective cover functor to obtain uniqueness
of $E_\infty$ structures. The first author thanks the Max-Planck
Institute and the mathematics department  in Bonn. The second
author was  partially supported by the \emph{Strategisk
Universitetsprogram i Ren Matematikk} (SUPREMA) of the Norwegian
Research Council. }

\maketitle

\begin{abstract}
We refine our earlier work on the existence and uniqueness of $E_\infty$
structures on $K$-theoretic spectra to show that at each prime $p$, the
connective Adams summand $\ell$ has a unique structure as
a commutative $\mathbb{S}$-algebra. For the $p$-completion $\ell_p$ we
show that the McClure-Staffeldt model for $\ell_p$ is equivalent as an
$E_\infty$ ring spectrum  to the connective cover of the periodic Adams
summand $L_p$. We establish a Bousfield equivalence between the
connective cover of the Lubin-Tate spectrum $E_n$ and $\bpn$.
\end{abstract}

\section*{Introduction}
The aim of this short note is to establish the uniqueness of $E_\infty$
structures on connective covers of certain periodic commutative
$\mathbb{S}$-algebras $E$, most prominently for the connective
$p$-complete Adams summand.  It is clear that the connective cover of an
$E_\infty$ ring spectrum inherits a $E_\infty$ structure; there is even a
 \emph{functorial} way of assigning a connective cover within the category
of $E_\infty$ ring spectra \cite[VII.3.2]{May}. But it is not obvious in
general that this $E_\infty$ multiplication is unique. 

Our main concern is with examples in the vicinity of $K$-theory; we apply
our uniqueness theorem to real and complex $K$-theory and their
localizations and completions and to the Adams summand and its completion. 

The existence and uniqueness of $E_\infty$ structures on the periodic
spectra $KU$, $KO$ and $L$ was established in \cite{AB&BR:HGamma} by means
of the obstruction theory for $E_\infty$ structures developed by
Goerss-Hopkins \cite{GH} and Robinson \cite{Rob}. Note however, that
obstruction theoretic methods would fail in the connective cases. Let $e$
be a commutative ring spectrum. If $e$ satisfies some K\"unneth and
universal coefficient properties \cite[proposition 5.4]{Rob}, then the
obstruction groups for $E_\infty$ multiplications consist of
Andr\'e-Quillen cohomology groups in the context of differential graded
$E_\infty$ algebras applied to the graded commutative $e_*$-algebra
$e_*e$. Besides problems with non-projectivity of $e_*e$ over $e_*$, the
algebra structures of $ku_*ku$, $ko_*ko$ and $\ell_*\ell$ are far from
being \'etale and therefore one would obtain non-trivial obstruction
groups. One would then have to identify actual obstruction classes in
these obstruction groups in order to establish the uniqueness of the given
$E_\infty$ structure -- but at the moment, this seems to be an intractable
problem. Thus an alternative approach is called for.

In Theorem~\ref{thm:Main} we prove that a unique $E_\infty$ structure
on $E$ gives rise to a unique structure on the connective cover if $E$
is obtained from some connective spectrum via a process of Bousfield
localization. In particular, we identify the $E_\infty$ structure on
the $p$-completed connective Adams summand $\ell_p$ provided by McClure
and Staffeldt in \cite{JM&RS} with the one that arises by taking the
unique $E_\infty$ structure on the periodic Adams summand $L = E(1)$
developed in~\cite{AB&BR:HGamma} and taking its connective cover.

Our Theorem applies as well to the connective covers of the Lubin-Tate
spectra $E_n$ and we prove in section~\ref{sec:bpn} that these spectra
are Bousfield equivalent to the truncated Brown-Peterson spectra
$\bpn$. Unlike other spectra that are Bousfield equivalent to $\bpn$,
such as the connective cover of the completed Johnson-Wilson spectrum,
$\widehat{E(n)}$, the connective cover of $E_n$ is calculationally
convenient. So far, only $\bp1 = \ell$ is known to have an $E_\infty$
structure, and we propose the connective cover of $E_n$ as an $E_\infty$
approximation of $\bpn$.

\section{$E_\infty$ structures on connective covers}
\label{sec:conncovers}

Let us first make explicit what we mean by uniqueness of
$E_\infty$ structures. We admit that this is an \emph{ad hoc} notion, but it
suffices for the examples we want to consider. 
\begin{defn}
In the following, we will say that an $E_\infty$ structure on some
homotopy commutative and associative ring spectrum $E$ is unique
if whenever there is a map of ring spectra $\varphi\: E'  \lra E$
from some other $E_\infty$ ring spectrum $E'$ to $E$ which induces
an isomorphism on homotopy groups, then there is a morphism in the
homotopy category of $E_\infty$ ring spectra $\varphi'\: E' \lra
E$ such that $\pi_*(\varphi) = \pi_*(\varphi')$.
\end{defn}

If $E$ and $F$ are spectra whose $E_\infty$ structure was provided
by the obstruction theory of Goerss and Hopkins \cite{GH}, then we can
compare our uniqueness notion with theirs. Note that examples of such
$E_\infty$ ring spectra include $E_n$ \cite[7.6]{GH}, $KO$, $KU$, $L$
and $\widehat{E(n)}$ \cite{AB&BR:HGamma}. In such cases the Hurewicz map
\begin{equation}\label{eqn:hur}
\Hom_{E_\infty}(E,F) \xrightarrow{h}
\Hom_{F_*\mathrm{-alg}}(F_*E,F_*)
\end{equation}
is an isomorphism. Assume that we have a mere ring map $\phi$ as
above between $E$ and $F$. This gives rise to a map of
$F_*$-algebras from $F_*E$ to $F_*$ by composing $F_*(\phi)$ with
the multiplication $\mu$ in $F_*F$. The left hand side in
\eqref{eqn:hur} denotes the derived space of $E_\infty$ maps from
$E$ to $F$. In presence of a universal coefficient theorem we have
$\Hom_{F_*\mathrm{-hom}}(F_*E,F_*) = [E,F]$, therefore the element
$\mu \circ F_*(\phi)$ gives rise to a homotopy class of maps of
ring spectra $\tilde{\phi}$ from $E$ to $F$. We can assume that we
have functorial cofibrant replacement $Q(-)$, hence we obtain a
ring map  $Q(\tilde{\phi})$ from $Q(E)$ to $Q(F)$. Via the
isomorphism \eqref{eqn:hur} this gives a map of $E_\infty$ spectra
from $Q(E)$ to $Q(F)$, $\Phi$, therefore we obtain a zigzag
$$
\xymatrix{
{Q(E)} \ar@{.>}[r]^{\Phi} \ar@{.>}[d]^{\sim} & {Q(F)} \ar@{.>}[d]^{\sim}\\
{E} \ar[r]^{\phi} & {F}}
$$
of weak equivalences of $E_\infty$ spectra from $E$ to $F$. Thus in
such cases our definition agrees with the uniqueness notion that is
natural in the Goerss-Hopkins setting.

For the rest of the paper we assume the following.
\begin{assump}
Let $E$ be a periodic commutative $\mathbb{S}$-algebra with periodicity
element $v \in E_*$ of positive degree. We will view $E$ as being obtained
from a connective commutative $\mathbb{S}$-algebra $e$ by Bousfield
localization at $e[v^{-1}]$ in the category of
$e$-modules. Furthermore we assume that the localization map induces
an isomorphism between the homotopy groups of $e$ and
the homotopy groups of the connective cover of $E$.
\end{assump}
Let us denote the connective cover functor from~\cite[VII.3.2]{May} by
$c(-)$. For any
$E_\infty$ ring spectrum $A$, there is a weakly equivalent commutative
$\mathbb{S}$-algebra $B(\mathbb{P},\mathbb{P}, \mathbb{L})(A)$, with
equivalence
$$
\lambda\: B(\mathbb{P},\mathbb{P}, \mathbb{L})(A) \xrightarrow{\simeq} A,
$$
in the $E_\infty$ category~\cite[XII.1.4]{EKMM}. Here
$B(\mathbb{P},\mathbb{P}, \mathbb{L})$ is a bar construction with respect
to the monad associated to the linear isometries operad $L$ and the monad
for commutative monoids in the category of $\mathbb{S}$-algebras 
$\mathbb{P}$. We
will denote the composite $B(\mathbb{P},\mathbb{P}, \mathbb{L}) \circ c$
by $\bar{c}$. For a commutative $\mathbb{S}$-algebra $R$ and an $R$-module
$M$, let $L^R_M(-)$ denote Bousfield localization at $M$ in the category
of $R$-modules and we denote the localization map by
$\sigma\: E \lra L_M^R(E)$ for any $R$-module $E$.

\begin{thm}\label{thm:Main}
Assume that we know that the $E_\infty$ structure on $E$ is
unique. Then 
the $E_\infty$ structure on $c(E)$ is unique.
\end{thm}

\begin{proof}
Each commutative $\mathbb{S}$-algebra can be viewed as an
$E_\infty$ ring spectrum. Let $e'$ be a model for the connective
cover $c(E)$, \ie, $e'$ is an $E_\infty$ ring spectrum with a map
of ring spectra $\varphi$ to $c(E)$, such that $\pi_*(\varphi)$ is
an isomorphism. Write $v \in e'_*$ for the isomorphic image of $v$
under the inverse of $\pi_*(\varphi)$. As $\varphi$ is a ring map
it will induce a ring map on the corresponding Bousfield
localizations. But as the $E_\infty$ structure on $E$ is unique by
assumption, this map can be replaced by an equivalent equivalence,
$\xi$, of $E_\infty$ ring spectra. We abbreviate
$B(\mathbb{P},\mathbb{P}, \mathbb{L})(e')$ to $B(e')$. We consider
the following diagram whose dotted lines provide a 
zigzag of $E_\infty$ equivalences and hence a map in the homotopy
category of $E_\infty$ ring spectra. 
\[
\xymatrix{
{} & {c(B(e'))} \ar@{.>}[dd]^(0.4){c(\sigma)}
\ar@{.>}[dl]_{\epsilon} & {} & {} & {c(e)} \ar@{.>}[dd]^(0.4){c(\sigma)}
\ar@{.>}[dl]_{\epsilon}\\
{B(e')}
\ar@{.>}[rr]^(0.4){\lambda} \ar[dd]^(0.4){\sigma} & {} & {e'} \ar@{<~>}[r] & {e}
\ar[dd]^(0.4){\sigma} & {} \\
{} & {c(L^{B(e')}_{B(e')[v^{-1}]}(B(e')))}
\ar@{.>}[rrr]^{c(\xi)} \ar[dl]_{\epsilon} & {} & {} & {c(E)}\ar[dl]_{\epsilon} \\
{L^{B(e')}_{B(e')[v^{-1}]}(B(e'))}
\ar[rrr]^{\xi}& {} & {} & {E} & {}
}
\]
\end{proof}

Real and complex $K$-theory, $ko$ and $ku$, have $E_\infty$
structures obtained using algebraic $K$-theory models~\cite[VIII,
\S2]{May}. The connective Adams summand $\ell$ has an $E_\infty$
structure because it is the connective cover of $E(1)$. In the
following we will refer to these models as the standard ones. The
$E_\infty$ structures on $KO$, $KU$ and $E(1)$ are unique
by~\cite[theorems~7.2,~6.2]{AB&BR:HGamma}. In all of these cases,
the periodic versions are obtained by Bousfield
localization~\cite[VIII.4.3]{EKMM}.

\begin{cor}\label{cor:ktheory}
The $E_\infty$ structures on $ko$, $ku$ and $\ell$ are unique. 
\end{cor}

In~\cite{JM&RS}, McClure and Staffeldt construct a model for the
$p$-completed connective Adams summand using algebraic $K$-theory
of fields. Let $\tilde{\ell} = K(\mathbf{k}')$, the algebraic
$K$-theory spectrum of $\mathbf{k}' = \bigcup_i \F_{q^{p^i}}$,
where $q$ is a prime which generates the $p$-adic units $\Z_p^\times$.
Then the $p$-completion of $\tilde{\ell}$ is additively equivalent
to the $p$-completed connective Adams summand $\ell_p$
\cite[proposition~9.2]{JM&RS}. For further details see
also~\cite[\S 1]{Ausoni}. Note that the $p$-completion $\ell_p$
inherits an $E_\infty$ structure from $\ell$ because $p$-completion is
Bousfield localization with respect to $H\F_p$ and therefore preserves
commutative $\mathbb{S}$-algebras \cite[VIII.2.2]{EKMM}. An  
\emph{a priori} different model for
the $p$-completion of the connective Adams summand can be obtained
by taking the connective cover of the $p$-complete periodic version
$L = E(1)$. This is consistent with the statement of
Corollary~\ref{cor:ktheory} because $p$-completion and Bousfield 
localization of $\ell$ in the category of $\ell$-modules with respect
to $L$ are compatible in the following sense. Consider
$\ell = \bar{c}(L)$ and its $p$-completion
$$
\lambda_\ell\: \ell \lra \ell_p = (\bar{c}(L))_p.
$$
The $p$-completion map $\lambda$ is functorial in the spectrum,
therefore the following diagram of solid arrows commutes.
$$
\xymatrix{
{\ell = \bar{c}(L)} \ar[rr]^{\lambda_\ell} \ar[dr] & &
{\ell_p = \bar{c}(L)_p} \ar[dr]
\ar@{.>}[rr] & & {\bar{c}(L_p)} \ar[dl] \\
 & {L} \ar[rr]^{\lambda_L} & & {L_p} &
}
$$
The universal property of the connective cover functor ensures
that there is a map in the homotopy category of commutative
$\mathbb{S}$-algebras from $\ell_p$ to $\bar{c}(L_p)$ which is
a weak equivalence. In the following we will not distinguish
$\ell_p$ from $\bar{c}(L_p)$ anymore and denote this model
simply by $\ell_p$.
\begin{prop}\label{prop:Comparison}
The McClure-Staffeldt model $\tilde{\ell}_p$ of the $p$-complete
connective Adams summand is equivalent as an $E_\infty$ ring
spectrum to $\ell_p$.
\end{prop}

\begin{rem}\label{prop:Eqt}
If $E$ is a commutative $\mathbb{S}$-algebra with naive $G$-action
for some group $G$, then neither the connective cover functor
$\bar{c}(-)$ nor Bousfield localization of $E$ has to commute with
taking homotopy fixed points. As an  example, consider connective
complex $K$-theory $ku$ with the conjugation action by $C_2$. The
homotopy fixed points $ku^{hC_2}$ are not equivalent to $ko$, but
on the periodic versions we obtain $KU^{hC_2} \simeq KO$.
\end{rem}

\begin{proof}[of Proposition~\ref{prop:Comparison}]
Consider the algebraic $K$-theory model for connective complex $K$-theory,
$ku = K(\mathbf{k})$, with $\mathbf{k} = \bigcup_i \F_{q^{p^i(p-1)}}$.
The canonical inclusions $\F_{q^{p^i}} \hookrightarrow \F_{q^{p^i(p-1)}}$
assemble into a map $j\:\mathbf{k}' \lra \mathbf{k}$. The Galois group
$C_{p-1}$ of $\mathbf{k}$ over $\mathbf{k}'$ acts on $\mathbf{k}$ and
induces an action on algebraic $K$-theory. As $\mathbf{k}'$ is fixed
under the action of $C_{p-1}$ there is a factorization of $K(j)_p$ as
$$
\xymatrix{
{K(\mathbf{k}')_p} \ar[dr]^{i} \ar[rr]^{K(j)_p} & & {K(\mathbf{k})_p} \\
& {K(\mathbf{k})_p^{hC_{p-1}} } \ar[ur] & }
$$
and $i$ yields a weak equivalence of commutative $\mathbb{S}$-algebras,
where $K(\mathbf{k})_p^{hC_{p-1}}$ is a model for the connective
$p$-complete Adams summand which is weakly equivalent to $\tilde{\ell}_p$
 (see \cite[\S 1]{Ausoni}).

Consider the composition of the following chain of maps between commutative
$\mathbb{S}$-algebras:
$$
K(\mathbf{k}')_p \xrightarrow{\;i\;} (K(\mathbf{k})_p)^{hC_{p-1}}
                                         \lra K(\mathbf{k})_p \lra KU_p.
$$
The target $KU_p$ is as well the target of the map $ \bar{c}(KU_p)
\lra KU_p$. Note that the universal property of $\bar{c}(-)$
yields a zigzag $\varsigma\: K(\mathbf{k})_p \leftrightsquigarrow
\bar{c}(KU_p)$ of equivalences between $K(\mathbf{k})_p$ and
$\bar{c}(KU_p)$ in the category of commutative
$\mathbb{S}$-algebras.

As $KU_p$ is the Bousfield localization of $K(\mathbf{k})_p$ in the
category of $K(\mathbf{k})_p$-modules with respect to the Bott element,
$$
KU_p = L^{K(\mathbf{k})_p}_{K(\mathbf{k})_p[\beta^{-1}]}K(\mathbf{k})_p,
$$
it inherits the $C_{p-1}$-action on $K(\mathbf{k})_p$. The functoriality
of the connective cover lifts this action to an action on $\bar{c}(KU_p)$.

The connective cover functor is in fact a functor in the category of
commutative $\mathbb{S}$-algebras with multiplicative naive $G$-action for 
any group
$G$. To see this we  have to show that the map $\bar{c}(A) \lra A$ is
$G$-equivariant
if $A$ is a commutative $\mathbb{S}$-algebra with an underlying 
naive $G$-spectrum. The functor $B(\mathbb{P},\mathbb{P}, \mathbb{L})$ does
not cause any problems. Proving the claim for the functor $c$ involves
chasing the definition given in~\cite[VII,~\S3]{May}.

The prespectrum underlying $c(A)$ applied to an inner product
space $V$ is defined as $T(A_0)(V)$, where $A_0$ is the zeroth space of the
spectrum $A$ and $T$ is a certain bar construction involving
suspensions and a monad consisting of the product of a fixed
$E_\infty$ operad with the partial operad of little convex bodies
$\mathcal{K}$. For a fixed $V$ the suspension $\Sigma^V$ and the
operadic term  $\mathcal{K}_V$ are used.  As the $G$-action is
compatible with the $E_\infty$ and the additive structure of $A$,
the evaluation map $T(A_0)(V) \lra A(V)$ is $G$-equivariant. For
varying $V$, these maps constitute a map of prespectra and its
adjoint on the level of spectra is $c(A) \lra A$. As the
spectrification functor preserves $G$-equivariance, the claim
follows. Therefore the resulting zigzag $\varsigma\:
K(\mathbf{k})_p \leftrightsquigarrow \bar{c}(KU_p)$ is
$C_{p-1}$-equivariant and we obtain an induced zigzag on homotopy
fixed points,
$$
\varsigma^{hC_{p-1}}\: (K(\mathbf{k})_p)^{hC_{p-1}}
                   \leftrightsquigarrow (\bar{c}(KU_p))^{hC_{p-1}}.
$$
As $\varsigma$ is an isomorphism in the homotopy category and is
$C_{p-1}$-equivariant, $\varsigma^{hC_{p-1}}$ yields an isomorphism
as well.
\end{proof}

\section{Connective Lubin-Tate spectra} \label{sec:bpn}

Goerss and Hopkins proved in \cite{GH} that the Lubin-Tate spectra
$E_n$ with
$$
{(E_n)}_* = W(\F_{p^n})[[u_1, \ldots, u_{n-1}]][u^{\pm 1}]
\quad
\text{with } |u_i| = 0 \text{ and } |u| = -2
$$
possess unique $E_\infty$ structures for all primes $p$ and all
$n \geq 1$. The connective cover $c(E_n)$ has coefficients
$$
{(c(E_n))}_* = W(\F_{p^n})[[u_1, \ldots, u_{n-1}]][u^{-1}]
\quad \text{ with } |u_i| = 0 \text{ and } |u| = -2.
$$
Of course $\bar{c}(E_n)[(u^{-1})^{-1}] \sim E_n$.

The spectra $\bpn$ can be built from the Brown-Peterson spectrum
$\BP$ by killing all generators of the form $v_m$ with $m>n$ in
$\BP_* = \Z_{(p)}[v_1,v_2,\ldots]$. Using for instance Angeltveit's
result~\cite[theorem 4.2]{A} one can prove that the $\bpn$ are
$A_\infty$ spectra and from~\cite{AB&AJ} it is known that this
$\mathbb{S}$-algebra structure can be improved to an $\MU$-algebra
structure. On the other hand, Strickland showed in~\cite{St} that
$\bpn$ with $n \geq 2$ is not a homotopy commutative $\MU$-ring
spectrum for $p=2$. We offer $c(E_n)$ as a replacement for the
$p$-completion $\bpn_p$ of $\bpn$.

We also need to recall that in the category of $\MU$-modules, $E(n)$
is the Bousfield localization of $\bpn$ with respect to $\bpn[v_n^{-1}]$,
hence by~\cite{EKMM} it inherits the structure of an $\MU$-algebra
and the natural map $\bpn \lra E(n)$ is a morphism of $\MU$-algebras.
Furthermore, the Bousfield localization of $E(n)$ with respect to
the $MU$-algebra $K(n)$ is the $I_n$-adic completion $\widehat{E(n)}$,
which was shown to be a commutative $\mathbb{S}$-algebra
in~\cite{AB&BR:HGamma}, and the natural map $\widehat{E(n)} \lra E_n$
is a morphism of commutative $\mathbb{S}$-algebras, see for
example~\cite[example~2.2.6]{AB&BR:BNG}. Thus there is a morphism of
ring spectra $\bpn \lra E_n$ which lifts to a map $\bpn \lra c(E_n)$.
\begin{prop}\label{prop:bpn->cEn}
The spectra $\bpn$ and $\bpn_p$ are Bousfield equivalent to  $c(E_n)$.
\end{prop}
\begin{proof}
On coefficients, we obtain a ring homomorphism
${(\bpn_p)}_* \lra {(c(E_n))}_*$ which on homotopy is given by
$$v_k \longmapsto
\begin{cases}
 u^{1-p^k}u_k & \text{ for } 1 \leq k \leq n-1,\\
 u^{1-p^n} & \text{ for } k = n. \end{cases}$$
extending the natural inclusion of the $p$-adic integers $\Z_p =
W(\F_{p})$  into  $W(\F_{p^n})$. This homomorphism is
induced by a map of ring spectra.

Recall from~\cite{AB:En-local} that $E(n)$ and $\widehat{E(n)}$
are Bousfield equivalent as $\mathbb{S}$-modules, and it follows
that $E_n$ is Bousfield equivalent to these since it is a finite
wedge of suspensions of $\widehat{E(n)}$.


If $X$ is a $p$-local spectrum with torsion free homotopy
groups then its $p$-completion $X_p$ is Bousfield equivalent
to $X$, \ie, $\langle X_p\rangle = \langle X\rangle$. This
follows using the cofibre triangles (in which $M(p)$ is the
mod~$p$ Moore spectrum and the circled arrow indicates a map of degree
one)
\[
\xymatrix{
X \ar[rr]^{p} && X \ar[dl] \\
& X\wedge M(p) \ar|\circ[ul] &
}
\quad
\xymatrix{
X_p \ar[rr]^{p} && X_p\ar[dl] \\
& X\wedge M(p) \ar|\circ[ul] &
}
\]
together with the fact that the rationalization $p^{-1}X$
is a retract of $p^{-1}(X_p)$. In particular, we have
$\langle \bpn_p\rangle = \langle \bpn \rangle$ and
$\langle E(n)_p\rangle = \langle E(n) \rangle$.

From~\cite[theorem~2.1]{Rav:loc}, the Bousfield class of $\bpn$ is
$$
\langle \bpn\rangle = \langle E(n)\rangle \vee \langle H\F_p\rangle.
$$
There is a cofibre triangle
\[
\xymatrix{
\Sigma^{2}c(E_n)\ar[rr]^{u^{-1}} && c(E_n)\ar[dl] \\
& {HW(\F_{p^n})[[u_1,\ldots,u_{n-1}]]}\ar|\circ[ul] &
}
\]
in which $HW(\F_{p^n})[[u_1,\ldots,u_{n-1}]]$ is  the Eilenberg-MacLane
spectrum. More generally we can construct a family of Eilenberg-MacLane
spectra with $W(\F_{p^n})[[u_1,\ldots,u_k]]$ as coefficients 
for $k=0,\ldots,n-1$ which are  related by cofibre triangles
\[
\xymatrix{
{\scriptstyle HW(\F_{p^n})[[u_1,\ldots,u_k]]}\ar[rr]^{u_k} && 
{\scriptstyle HW(\F_{p^n})[[u_1,\ldots,u_k]]}\ar[dl] \\
& {\scriptstyle HW(\F_{p^n})[[u_1,\ldots,u_{k-1}]]}\ar|\circ[ul] &
}
\]
such that for $k=0$ we obtain $HW(\F_{p^n})$. With the help of these cofibre sequences
we can identify  $ \langle HW(\F_{p^n})[[u_1,\ldots,u_k]] \rangle$ with 
$\langle HW(\F_{p^n})[[u_1,\ldots,u_{k-1}]] \rangle  \vee
\langle HW(\F_{p^n})[[u_1,\ldots,u_k]][u_k^{-1}] \rangle$.

In general, if $R$ is a commutative ring, then the ring of finite
tailed Laurent series $R((x))$ is faithfully
flat over $R$ and  therefore we have
$$ \langle HR((x)) \rangle = \langle HR \rangle.$$
Using this auxiliary fact we inductively get that
$$ \langle HW(\F_{p^n})[[u_1,\ldots,u_k]] \rangle = \langle
HW(\F_{p^n})[[u_1,\ldots,u_{k-1}]] \rangle.$$
This reduces the Bousfield class of $c(E_n)$ to $\langle E_n \rangle
\vee \langle HW(\F_{p^n})\rangle$. As $W(\F_{p^n})$ is a finitely generated
free $\Z_{p}$-module and as $\langle H\Z_{p}\rangle = \langle
H\Q \rangle \vee\langle H\F_p \rangle$ this leads to
\begin{align*} \langle c(E_n) \rangle & = \langle E(n) \vee H\Q \vee
  H\F_p \rangle \\
& = \langle E(n) \vee
  H\F_p \rangle = \langle BP\langle n\rangle  \rangle. \qedhere
\end{align*}
\end{proof}

\end{document}